\documentclass{amsart}
\usepackage{mathrsfs}
\usepackage{enumerate}

\newtheorem{theorem}{Theorem}

\newtheorem{proposition}{Proposition}[section]

\theoremstyle{definition}
\newtheorem{definition}{Definition}[section]

\theoremstyle{remark}
\newtheorem{remark}{Remark}[section]

\numberwithin{equation}{section}

\def\d{\text{d}}
\def\e{\text{e}}
\def\i{\text{i}}
\def\m{{\mathbf m}}
\def\x{{\mathbf x}}
\def\y{{\mathbf y}}
\def\l{{\mathbf l}}
\def\r{{\mathbf r}}

\begin{document}

\title[Diffusion processes and continued fractions]{Excursions of diffusion processes
and continued fractions}

\author{Alain Comtet}
\address{Laboratoire de Physique Th\'eorique et Mod\`eles Statistiques, Universit\'e Paris-Sud, B\^atiment 100, 91405 Orsay Cedex, France \\
              Universit\'e Pierre et Marie Curie--Paris 6, 4 Place Jussieu, 75252 Paris Cedex 5, France} 
              \email{alain.comtet@u-psud.fr}        %  \\
\author{Yves Tourigny}
\address{School of Mathematics\\
        University of Bristol\\
        Bristol BS8 1TW, United Kingdom}
\email{y.tourigny@bristol.ac.uk}

\thanks{We thank Professor Marc Yor for  many useful
discussions, and for his comments on earlier versions of this work. The second author also acknowledges the support
of the Engineering and Physical Sciences Research Council
(United Kingdom) under Grant GR/S87461/01.}

%    General info

\date{\today}


\begin{abstract}
It is well-known that the excursions of a one-dimensional diffusion process
can be studied by considering a certain Riccati equation associated with the process. We show that,
in many cases of interest, the Riccati equation can be solved in terms of an infinite continued fraction.
We examine the
probabilistic significance of the expansion. To illustrate our results, we discuss some examples of diffusions
in deterministic and in random environments. 
\end{abstract} 

\maketitle

\section{Introduction}
\label{intro}
Recently, Marklof {\em et al.} \cite{MTW1,MTW2} studied the random continued fraction
\begin{equation}
\label{stieltjesContinuedFraction}
\cfrac{2 \lambda}{u_1  + \cfrac{2 \lambda}{u_2+ \cfrac{2 \lambda}{u_3 + \cdots}}}\,, \quad \lambda \in {\mathbb C}\,,
\end{equation}
where 
the $u_n$ are independent random variables with the same gamma distribution, i.e. for every Lebesgue-measurable
set $A \subset {\mathbb R}_+$,
\begin{equation}
{\mathbb P} \left ( u_n \in A \right ) = \int_A \frac{1}{2^{\mu} \Gamma (\mu)} y^{\mu-1} \e^{-y/2}\,\d y\,, \quad \mu > 0\,.
\label{gammaDistribution}
\end{equation}
An earlier work of Letac \& Seshadri \cite{LS} had shown that, in the case
of positive $\lambda$, the continued fraction is a generalised inverse Gaussian random variable. 
One of the results obtained by Marklof {\em et al.} is an explicit formula for the probability density function  for every
complex $\lambda$.

The genesis of the work reported here was the curious observation that the {\em same distribution}
--- not the continued fraction itself---
appears also in the articles
of Bouchaud {\em et al.} \cite{BCGL} ($\lambda <0$) and Kawazu \& Tanaka \cite{KT} ($\lambda > 0$)
on diffusion in a Brownian environment with drift. More precisely, consider a process $X$ with 
infinitesimal
generator
\begin{equation}
{\mathscr G}_W := \frac{a}{2}  \e^{-2 W} \frac{\d}{\d x} \left [ \e^{2 W}  \frac{\d}{\d x} \right ]\,,
\label{randomGenerator}
\end{equation} 
where
$$
W(x) := \int \frac{b}{a}\,\d x
$$
and the functions $a$ and $b$ are the instantaneous variance and instantaneous drift respectively. 
Many fundamental quantities associated with the process can be expressed in terms of the solutions of the equation
\begin{equation}
{\mathscr G}_W \phi(\cdot,\lambda) = \lambda \phi(\cdot,\lambda)\,.
\label{diffusionEquation}
\end{equation}
The order of the equation may be lowered by introducing the Riccati variable
$$
U := \frac{\phi'}{\phi}\,,
$$
where the prime symbol denotes differentiation with respect to $x$.
This yields the {\em Riccati equation}
\begin{equation}
\frac{\d U}{\d x} + U^2 + 2 W'(x) \,U = 2 \lambda/a(x)\,.
\label{riccatiEquation}
\end{equation}
When the state space is ${\mathbb R}$ and $W$ is Brownian motion with drift, 
this Riccati
equation (understood in the sense of Stratonovich \cite{RY})
admits a positive stationary solution whose law is the same as that of the continued fraction (\ref{stieltjesContinuedFraction}).
This formulation in terms of the Riccati variable goes back to the pioneering work of Frisch \& Lloyd \cite{FL} and has since been used
extensively, both in the physics \cite{BCGL,MC} and the mathematics \cite{Ko,KT} literatures.

In this work, we show that, for many diffusions--- be they in a fixed or a random environment $W$--- a solution of this Riccati equation
may be found in terms of the infinite continued fraction
\begin{equation}
U (x,\lambda) = u_0 (x) + \cfrac{2\lambda/a(x)}{u_1(x) + \cfrac{2 \lambda/a(x)}{u_2(x)+ \cdots}}\,.
\label{continuedFraction}
\end{equation}

Continued fractions arise in the study of many random processes on a discrete state space, notably birth-and-death processes (see
for instance \cite{BR,DFPS,FG,GP,KM} and the references therein) and random
walks in a random environment \cite{APS, BAO, De};
their occurence in the context of processes on a continuous state space is comparatively rare, though not unknown \cite{BY}.
In the remainder of this introduction, we provide a brief explanation of the origin and significance of the continued fraction
in the context of diffusion processes, and summarise our main results.

\subsection{Probabilistic interpretation of the Riccati variable}
\label{riccatiSubsection}
The equation (\ref{diffusionEquation}), supplemented with the boundary conditions that may be required in order to specify the process $X$ uniquely,
has two non-negative solutions, say $\phi_{-} (\cdot,\lambda)$ and $\phi_+ (\cdot,\lambda)$, that are of particular importance:  they are characterised (up to a constant factor) by the fact that $\phi_-(\cdot,\lambda)$ is {\em non-decreasing} and $\phi_+(\cdot,\lambda)$ {\em non-increasing}.  
Let $x$ denote the {\em starting point} of the diffusion and let $y$ be any other point in the state space. Then the {\em first hitting time} $H(y)$ of $y$
is the random variable
defined by 
\begin{equation}
H(y) := \inf \left \{ t \ge 0 \,:\; X_t = y \right \}
\label{hittingTime}
\end{equation}
and we have the following well-known formula for its Laplace transform \cite{BS}:
\begin{equation}
{\mathbb E}_x \left ( \e^{-\lambda H(y)} \right ) = \begin{cases}
\frac{\phi_-(x,\lambda)}{\phi_-(y,\lambda)} & \text{if $x \le y$} \\
 & \\
\frac{\phi_+(x,\lambda)}{\phi_+(y,\lambda)} & \text{if $x \ge y$} 
\end{cases}\,.
\label{laplaceTransformOfHittingTime}
\end{equation}
We deduce
\begin{equation}
U_{\pm} (x,\lambda) := \frac{\phi_{\pm}'(x,\lambda)}{\phi_{\pm}(x,\lambda)} = \frac{\d}{\d y} {\mathbb E}_y \left ( \e^{-\lambda H(x)} \right ) \Bigl |_{y=x \pm} \,.
\label{pitmanYorFormula1}
\end{equation}
This expresses two particular solutions of the Riccati equation--- one positive ($U_-$) and one negative ($U_+$)--- in simple probabilistic terms. 
The choice of sign in this notation is somewhat disconcerting,
but its justification will soon become manifest.

Following Pitman \& Yor \cite{PY}, one can gain further insight by
defining a local time process $L(x)$ at $x$ by
\begin{equation}
L_t (x) := \lim_{\varepsilon \rightarrow 0+} \frac{a(x)}{2 \varepsilon} \int_0^t {\mathbf 1}_{(x-\varepsilon,x+\varepsilon)} \left ( X_{\tau} \right )\,\d \tau \,.
\label{localTime}
\end{equation}
This process is proportional to the time that $X$ spends in the vicinity of the starting point up to time $t$. Let $\zeta$ denote the lifetime
of $X$. If $X$ is recurrent (respectively transient), then $L_{\zeta}(x)$ is infinite (respectively finite) almost surely.
The {\em inverse local time} $L^{-1}(x)$ at
the starting point $x$ is the process defined by
\begin{equation}
L_t^{-1} (x) := \inf \left \{ \tau \,:\; L_{\tau}(x) > t \right \}\,,\;\; 0 \le t < L_{\zeta}(x)\,.
\label{inverseLocalTime}
\end{equation}
It is a measure of the time that elapses before $X$ has spent a total time $t$ in the vicinity of its starting point. 
Using the trivial identity
$$
t = \int_0^t {\mathbf 1}_{(-\infty,x]} \left ( X_{\tau} \right )\,\d \tau + \int_0^t {\mathbf 1}_{(x,\infty)} \left ( X_{\tau} \right )\,\d \tau\,,
$$
we can write
\begin{equation}
L^{-1} (x) = T_-(x) + T_{+}(x)
\label{occupationFormula}
\end{equation}
where
\begin{equation}
\notag
T_{-,t} (x) = \int_0^{L_t^{-1} (x)} {\mathbf 1}_{(-\infty,x]} \left ( X_{\tau} \right )\,\d \tau \;\;\text{and}\;\;
T_{+,t} (x) = \int_0^{L_t^{-1} (x)} {\mathbf 1}_{(x,\infty)} \left ( X_{\tau} \right )\,\d \tau\,.
\end{equation}
The processes $T_{-}(x)$ and $T_{+}(x)$ have obvious interpretations as occupation times;
they are subordinators, i.e. non-decreasing L\'{e}vy processes, and we have
\begin{equation}
{\mathbb E} \left ( \exp \left [ -\lambda T_{\pm,t}(x) \right ] \right ) = \exp \left [ - t \psi_{\pm} (x,\lambda) \right ]\,, \quad \lambda > 0\,,
\label{laplaceExponent}
\end{equation}
for some function $\psi_{\pm}(x,\cdot)$ called the {\em Laplace exponent} of $T_{\pm} (x)$ \cite{Be1}. By the L\'{e}vy--Khintchine formula, the Laplace
exponent can be written in the form
\begin{equation}
\psi_{\pm}(x,\lambda) = \int_{(0,\infty]} \left ( 1- \e^{-\lambda y} \right ) \,\nu_{\pm}(x,\d y) 
\label{levyKhintchineFormula}
\end{equation}
for some measure $\nu_{\pm}(x)$ on $(0,\infty]$ such that
$$
\int_{(0,\infty)} \min \{1,y\}\,\nu_{\pm} (x,\d y) < \infty\,.
$$

The measure $\nu_{\pm}(x,\cdot)$ is called the {\em L\'{e}vy measure} of the process $T_{\pm}(x)$; it describes the distribution of the heights
of its jumps. 
Each jump corresponds to an {\em excursion} of the process $X$, i.e. to 
the path followed by the process between two successive visits to the starting point $x$, and the jump height is the duration of the
excursion. The $+$ (respectively $-$) case picks out the ``upward'' (respectively ``downward'') excursions, i.e. those whose path lies entirely above (respectively below or at) $x$. In the transient case, one or both of $\nu_+(x,\cdot)$ and $\nu_-(x,\cdot)$ will have an atom at infinity, reflecting the fact that $X$ begins an
excursion of infinite duration at the time of its last visit to $x$. Pitman \& Yor \cite{PY} show that
\begin{equation}
\psi_{\pm}(x,\lambda) = \mp \frac{1}{2} U_{\pm}(x,\lambda) \,.
\label{pitmanYorFormula}
\end{equation}
Knowing $U_{\pm}$, one can therefore recover the L\'{e}vy measure $\nu_{\pm}(x,\cdot)$ by inverting a Laplace transform.

\subsection{ Stieltjes functions}
\label{stieltjesSubsection}

Knight \cite{Kn1} and Kotani \& Watanabe \cite{KW} make the observation that
Krein's theory of strings implies the existence of a
measure, say $\sigma_{\pm}(x)$, on $[0,\infty)$ such that
\begin{equation}
U_{\pm} (x,\lambda) 
= \mp 2 \lambda \int_{[0,\infty)} \frac{\sigma_{\pm}(x,\d z)}{\lambda+z}\,.
\label{knightFormula}
\end{equation}
So, knowing $U_{\pm}$, we can obtain $\sigma_{\pm}$ (and hence $\nu_{\pm}$) from
\begin{equation}
\sigma_{\pm} (x,\{0\}) = \frac{\mp 1}{2} U_{\pm} (x,0)
\label{inversionFormula1}
\end{equation}
and the {\em Stieltjes--Perron inversion formula}
\begin{equation}
\sigma_{\pm} (x,A) = \lim_{\varepsilon \rightarrow 0+} \frac{1}{2 \pi} \int_A \text{Im} \left [ \frac{\mp U_{\pm}(x,-\lambda-\i \varepsilon)}{-\lambda-\i \varepsilon} \right ]\,\d \lambda
\label{inversionFormula2}
\end{equation}
which holds for every $\sigma_\pm$-measurable set $A \subset (0,\,\infty)$.
Then 
\begin{equation}
\nu_{\pm} \left ( x,\{\infty\} \right ) = \sigma_{\pm} \left ( x,\{0\} \right ) \;\;\text{and}\;\;
\nu_{\pm} (x,\d y) = \int_0^\infty z \e^{-y z} \sigma_{\pm} (x,\d z)\,\d y\,.
\label{stieltjesMeasure}
\end{equation}

A function of the form
$$
S(\omega) = \int_0^\infty \frac{\sigma(\d z)}{\omega+z}\,, \quad \omega \in {\mathbb C} \backslash {\mathbb R_-}\,,
$$
is called a {\em Stieltjes function} (or {\em transform}).
Stieltjes \cite{St} considered the problem of recovering a measure on ${\mathbb R}_+$ from its moments, and 
the continued fraction
\begin{equation}
\cfrac{1}{m_0 \omega + \cfrac{1}{\ell_1+\cfrac{1}{m_1 \omega + \cfrac{1}{\ell_2 + \cdots}}}}
\label{stieltjesExpansion}
\end{equation}
played an important part in his solution of this moment problem, as we proceed to explain.
Suppose 
that the measure $\sigma$ satisfies the following moment condition:
\begin{equation}
\text{For every non-negative integer $n$},\;\; \int_0^\infty z^{n} \sigma (\d z)  \;\;\text{exists}\,.
\label{stieltjesMomentCondition}
\end{equation}
For simplicity, suppose also that $\sigma$ has infinitely many points of growth. Then $S$ has an asymptotic expansion
in decreasing powers of $\omega$, namely
$$
S(\omega) \sim \sum_{n=0}^\infty \left [ \int_0^{\infty} (-z)^{n} \sigma (\d z) \right ] \omega^{-n-1} \quad \text{as $\omega \rightarrow +\infty$}\,,
$$
and one may associate with $S$ two sequences $\{m_n\}$ and $\{\ell_n\}$ of positive numbers by requiring that
the finite truncations of (\ref{stieltjesExpansion}) have asymptotic expansions that agree with the above series up to some order \cite{Ak,BO}.
Importantly,  only the first $2n+1$ of the moments are required in order to calculate the coefficients $m_n$ and $\ell_n$.
Stieltjes showed how, knowing these coefficients, 
a measure with the required moments can be constructed \cite{Ak,NS,St}. Furthermore, he showed that if the series
$$
\sum_{n=1}^\infty m_n \quad \text{or} \quad \sum_{n=0}^\infty \ell_n
$$
diverges, then there is only one measure $\sigma$ with the given moments, and the infinite continued fraction converges to $S(\omega)$.

\subsection{Statement of the main results and outline of the paper}
\label{outlineSubsection}
Our first task is to bring out the relationship between 
the continued fraction of Stieltjes and our own.  To this end, and with a view to making the paper reasonably self-contained, we recall
in \S \ref{kreinSection} some of the results from Krein's theory of strings that are most relevant to diffusion processes.
In Krein's terminology,
Equation (\ref{stieltjesExpansion}) is the characteristic function of a string with a discrete
distribution of masses $m_n$ such that $m_0$ is at $0$ and the spacing between $m_{n}$ and $m_{n-1}$ is $\ell_n$.
Then, in \S \ref{proofOfTheorem1Section}, we prove the 
\begin{theorem}
Let $X$ be a non-singular diffusion process started at $x$. Suppose that $\sigma_{\pm}(x,\cdot)$ has infinitely
many points of growth and that the following moment condition holds:
\begin{equation}
\tag{M}
\text{For every positive integer $n$,} \quad 
\int_0^{\infty} z^{-n} \sigma_{\pm} (x,\text{\em d} z)  < \infty\,.
\label{momentCondition}
\end{equation}
Then there
is a sequence of positive numbers 
$$
\mp u_{1,\pm}(x), \,\mp u_{2,\pm}(x),\, \ldots
$$
such that, for every, $n=0,\,1,\,2,\,\ldots$
$$
U_{\pm}(x,\lambda) -U_{\pm}(x,0)= \cfrac{2 \lambda/a(x)}{u_{1,\pm}(x)+\cfrac{2 \lambda/a(x)}{u_{2,\pm}(x) + \cdots+\cfrac{2 \lambda/a(x)}{u_{n,\pm}(x)}}}
\;+ O(\lambda^{n+1})
$$
as $\lambda \rightarrow 0+$.

Conversely, suppose that
the identity
$$
U_{\pm}(x,\lambda) - U_{\pm}(x,0) = \cfrac{2 \lambda/a(x)}{u_{1,\pm} (x) + \cfrac{2 \lambda/a(x)}{u_{2,\pm}(x)+ \cdots}}
$$
holds for some infinite sequence $\mp u_{1,\pm}(x)$, $\mp u_{2,\pm}(x), \,\ldots$ of positive numbers such that the series
$$
\sum_{n=1}^\infty u_{n,\pm}(x)
$$
diverges. Then $\sigma_\pm(x,\cdot)$ has infinitely many points of growth and the moment condition (\ref{momentCondition}) holds. 
\label{stieltjesTheorem}
\end{theorem} 

By using (\ref{stieltjesMeasure}) and Fubini's
theorem, it is easy to see that the moment condition
(\ref{momentCondition}) can be written in the equivalent form 
\begin{equation}
\text{For every positive integer $n$,} \quad 
\int_0^{\infty} y^{n} \nu_{\pm} (x,\text{\em d} y)  < \infty\,.
\label{alternativeMomentCondition}
\end{equation}
So when the coefficients are positive,
we immediately glean some information about the tail behaviour of the L\'{e}vy measure.
In this case, the continued fraction (\ref{continuedFraction}) is related to
that used by Stieltjes in his study of the moment problem, and
the Riccati solution $U_{\pm}(x,\lambda)$ has an expansion in increasing powers of $\lambda$, from which
the continued fraction coefficients can be computed by a well-known algorithm \cite{BO}.

The moment condition (\ref{momentCondition}), or its equivalent form (\ref{alternativeMomentCondition}), is by no means necessary
for the Riccati solution $U_{\pm}$ to have a continued fraction expansion {\em if we allow the coefficients $u_{n,\pm}$ to be of
arbitrary sign}. In such cases, $U_{\pm}(x,\lambda)$ is no longer asymptotic to a power series in $\lambda$.
Nevertheless, the continued fraction coefficients may be computed by an
algorithm, presented in a more general form in Common \& Roberts \cite{CR},  which we describe in  \S \ref{expansionSection}.

\begin{theorem}
Suppose that the Riccati equation (\ref{riccatiEquation}) has a solution expressible in the
form (\ref{continuedFraction}). Then, for $n \in {\mathbb N}$, $u_n$ solves the {\em homogeneous} Riccati equation
\begin{equation}
\frac{\text{\em d} u_n}{\text{\em d} x} + u_n^2  + 2 W_n'(x)\,  u_n = 0\,,
\label{homogeneousRiccatiEquation}
\end{equation}
where
\begin{equation}
W_{n}' = \frac{a'}{2a} - u_{n-1} - W_{n-1}'\,,\;\; W_0' := W'\,.
\label{Wrecurrence}
\end{equation}
\label{expansionTheorem}
\end{theorem}

In some cases, this expansion algorithm admits a probabilistic interpretation in terms of a sequence of diffusion processes, where the $n$th process
corresponds to the environment $W_n$. In \S \ref{probabilisticSection}, we study the relationship between adjacent processes
in this sequence. This reveals a surprising connection with the Ciesielski--Taylor theorem generalised by Biane \cite{Bi,CT}, which
exhibits a large class of ordered pairs of diffusions such that some occupation time of the first process has the same
law as some hitting time of the second process.

In \S \ref{exampleSection}, we consider some simple examples of diffusions in a deterministic environment which illustrate these results. 

Finally, in \S \ref{randomSection}, we
return to
the problem that provided the initial motivation for this work, and consider the case where the environment $W$ is itself a diffusion process:
\begin{theorem}
Let $W$ be a standard Brownian motion with positive drift. Then the system of Stratonovich equations (\ref{homogeneousRiccatiEquation}-\ref{Wrecurrence})
has a stationary solution such that $u_0 = 0$ and, for $n \ge 1$, the $u_n$ are independent random variables with the same gamma distribution.
\label{randomTheorem}
\end{theorem}
Thus, we recover the continued fraction
of \cite{LS,MTW1,MTW2}.

As this outline indicates, the present paper deals exclusively with what may be called the {\em direct problem}: given the characteristics $a$ and $b$ of the diffusion,
together with appropriate boundary conditions,
we use a continued fraction expansion to compute the L\'{e}vy measure of the excursions. The {\em inverse problem}, namely that
of finding the diffusion given the L\'{e}vy measure, is also of great interest \cite{DY}. The relevance of our
results to the solution of this inverse problem will be developed in a separate publication.

\section{Krein's theory of strings and diffusions}
\label{kreinSection}

This section provides a succinct review of some well-known results concerning strings and diffusions.
Kotani and Watanabe's account of these topics in \cite{KW} is particularly well-suited to our purpose, and we follow them very closely.
The reader familiar with this material need only take note of Definition \ref{kreinDualDefinition} before proceeding to the next section.

\subsection{Strings}
\label{stringSubsection}
A {\em string}, say $\m$, is a function from $[0,\infty]$ to $[0,\infty]$ that is non-decreasing, right-continuous and
infinite at infinity.
Let $\m$ be a string  and suppose that $\m$ is not identically infinite. Set
$$
c := \inf \left \{ \x \,:\; \m(\x) > 0 \right \} \;\;\text{and}\;\; \ell := \sup \left \{ \x \,:\; \m(\x) < \infty \right \}\,.
$$
The number $\ell$ is called the {\em length} of the string. 

By setting $\m(0-) = 0$, we obtain from $\m$ a (Stieltjes) measure $\d \m$ on $[0,\infty)$.
Let $\omega > 0$ and denote by $\xi(\cdot,\omega)$ and
$\eta(\cdot,\omega)$ the (unique!) solutions of the following integral equations on  $[0,\ell)$:
$$
\xi (\x,\omega) = 1+\omega \int_{0}^{\x} \left [ \int_{0-}^{{y}+} \xi ({z},\omega) \, \d \m ({z}) \right ] \d {y} 
$$
and
$$
\eta (\x,\omega) = \x+\omega \int_{0}^{\x} \left [ \int_{0-}^{{y}+} \eta ({z},\omega) \, \d \m ({z}) \right ] \d {y} \,.
$$
$\xi$ (respectively $\eta$) can be viewed as the particular solution of the generalised
differential equation
\begin{equation}
\frac{\d^2}{\d \m \d \x}  f(\cdot,\omega) = \omega f (\cdot,\omega)\,,\;\; 0 < \x < \ell\,,
\label{kreinEquation}
\end{equation}
subject to the condition
$$
f (0,\omega) = 1\; \text{and} \; \frac{\d f}{\d \x} (0,\omega) = \m (0+) \,\omega \;\; \left (\text{respectively}\;f(0,\omega) = 0 \; \text{and}\;\frac{\d f}{\d \x} (0,\omega) = 1 \right )\,.
$$
We shall refer to $(\xi,\eta)$ as the pair of {\em fundamental solutions} associated with the string. 

The function
$$
{S} (\omega) := \lim_{\x \rightarrow \ell-} \frac{\eta(\x,\omega)}{\xi(\x,\omega)}
$$
is called the  {\em characteristic function} of the string.
It may be shown that ${S}$ is necessarily of the form
\begin{equation}
{S}(\omega) = c + \int_{[0,\infty)} \frac{\sigma (\d z)}{\omega+z}
\label{characteristicFunction}
\end{equation}
for some measure $\sigma$ on $[0,\infty)$ such that
$$
\int_{[0,\infty)} \frac{\sigma (\d y)}{1+y} < \infty\,.
$$
Conversely, Krein showed that every function of the form (\ref{characteristicFunction}) is the characteristic function of a unique string. 

The right-continuous inverse, denoted $\m^\ast$, of a string $\m$ is called the {\em dual string}
of $\m$. 
Denote by 
$$
S^{\ast} (\omega) = c^\ast + \int_{[0,\infty)} \frac{\sigma^{\ast} (\d z)}{z+\omega}
$$
the characteristic function of the dual string. Then
$$
c^{\ast} = \m (0+) \;\;\text{and}\;\;{S}^\ast (\omega)=\frac{1}{\omega {S}(\omega)}\,.
$$

\subsection{Generalised diffusion processes}
\label{gapSubsection}

Given a pair $\left (\m_-,\m_+ \right )$ of strings such that
$\m_-(0+) = 0$ and $\ell_{\pm} > 0$, let 
$\d \tilde{\m}_-$ be the image measure of $\d \m_-$ under the map $\x \mapsto -\x$ and set
$$
\d \m := \begin{cases}
\d \tilde{\m}_- & \text{on $(-\infty,0)$} \\
\d \m_+ & \text{on $[0,\infty)$}
\end{cases}\,.
$$
The support of $\d \m$, denoted $I_{\m}$, is  a subset of $(-\ell_-,\ell_+)$. Let $B$ be a standard
Brownian motion, denote by $L_t^B(\x)$ its local time and set
$$
\varphi(t) := \int_{\mathbb R} L_t^B (\x) \,\d \m\,.
$$
Then 
$$
{\mathbf X}_t := B_{\varphi^{-1}(t)}
$$
defines a Markov process on $I_{\m}$ whose lifetime is the first hitting time of $-\ell_-$ or $\ell_+$.
${\mathbf X}$ is called the {\em generalised diffusion} corresponding to the pair  $\left (\m_-,\m_+ \right )$.

\begin{definition}
A {\em one-sided} diffusion is a generalised diffusion such that either $\m_-$ or $\m_+$ is the zero string.
\label{oneSidedDefinition}
\end{definition}
Thus a generalised diffusion may be thought of as an ordered pair of one-sided diffusions.

The following non-standard definition will also be helpful later on:
\begin{definition}
The {\em Krein dual}  ${\mathbf X}^\ast$ of the generalised diffusion process
${\mathbf X}$ corresponding to the pair $\left (\m_-,\m_+ \right )$ is the generalised diffusion process corresponding to the
pair $\left (\m_-^\ast,\m_+^\ast \right )$.
\label{kreinDualDefinition}
\end{definition}

\subsection{Non-singular diffusion processes}
\label{diffusionSubsection}
Next, consider 
a diffusion process $X$ with generator (\ref{randomGenerator}) whose state space $I$ is an interval with left endpoint $l$ and
right endpoint $r$. For the sake of greater clarity, we shall, in this subsection, depart from our usual notation and
use $x_0$ instead of $x$ to denote the starting point of the process;
this leaves us free to use $x$ to denote some generic point in the closure of $I$. 
We say that $X$ is {\em non-singular} \cite{Kn2,KW} or {\em regular} \cite{BS} if,
for every $x_0,\,y \in I$, $H(y)$ is finite with positive probability. We shall show that
such a process is also a generalised diffusion process in the sense
of the previous subsection.

Associated with the diffusion is a {\em scale function}
$s$ and a {\em speed measure} $m$ given respectively by
$$
s(x) = \int_{x_0}^x \e^{-2 W (y)}\,\d y \;\;\text{and}\;\; m(\d x) = \frac{2}{a} \e^{2 W}\,\d x\,. 
$$
The infinitesimal generator may then by expressed in the equivalent form
$$
{\mathscr G}_W = \frac{\d}{\d m} \left ( \frac{1}{s'(x)} \frac{\d}{\d x} \right )\,, \;\; l < x < r\,. 
$$
In general,
the scale function and the speed measure are not enough to determine completely (up to a constant factor)
the non-negative monotonic solutions $\phi_-$ and $\phi_+$ of Equation (\ref{diffusionEquation}) that characterise the process; some boundary conditions must also be imposed.
Following the accepted terminology, we say that $r$ (respectively $l$) is {\em entrance}  if
$$
 \int_{x_0}^r \left [ s(x)-s(x_0) \right ] \, m(\d x)   \; \left ( \text{respectively}\;\int_l^{x_0} \left [ s(x_0)-s(x) \right ] \, m(\d x) \right ) < \infty
$$
and that $r$ (respectively $l$) is {\em exit}  if
$$
\int_{x_0}^r  m((x_0,x))\, s'(x)\,\d x  \; \left ( \text{respectively}\; \int_l^{x_0}  m((x,x_0)) s'(x)\,\d x \right ) < \infty\,.
$$
A point is called {\em non-singular} if it is both entrance and exit, and {\em natural} if it is neither.
A boundary
condition for $\phi_+(\cdot,\lambda)$ (respectively $\phi_-(\cdot,\lambda)$) is required only if $r$ (respectively $l$) is non-singular, and then only at $r$ (respectively $l$):
\begin{enumerate}
\item If $r$ (respectively $l$) $\in I$, then this Feller-type condition is
\begin{equation*}
\alpha_+ \phi_+ (r,\lambda) + \beta_+ s'(r)\,\phi_+'(r,\lambda) + \gamma_+ \lambda \phi_+(r,\lambda) = 0\,,
%\label{leftBoundaryCondition}
\end{equation*}
\begin{equation*}
\left ( \text{respectively} \;\;\alpha_-\phi_- (l,\lambda) + \beta_- s'(l)\,\phi_-'(l,\lambda) + \gamma_- \lambda \phi_-(l,\lambda) = 0\,, \right )
%\label{rightBoundaryCondition}
\end{equation*}
where 
$$
\alpha_\pm \ge 0\,,\;\beta_\pm > 0\,, \; \gamma_\pm \ge 0\;\;\text{and}\;\;\alpha_\pm + \beta_\pm + \gamma_\pm = 1\,.
$$

\item If $r$ (respectively $l$) $\notin I$ then
the killing condition is
$$
\phi_+(r-,\lambda) \; \left ( \text{respectively}\; \phi_- (l+,\lambda) \right ) = 0\,.
$$
\end{enumerate}

Next, we explain how to construct a pair of strings associated with $X$. Set $\x = s(x)$.
Then 
$$
{\mathscr G}_W= \frac{\d^2}{\d \m \d \x}\,,\;\; {\tt l} < \x < {\tt r}\,,
$$
where
$$
\m (\d \x) = m(\d x)\,, \;\; \l := \int_{x_0}^l s'(y)\,\d y\;\;\text{and}\;\;\r := \int_{x_0}^r s'(y)\,\d y\,.
$$
Let $\m_\pm$ be as follows: for $0 \le \x < \r$ (respectively $0 \le \x < -\l$), set
$$
\m_+ (\x) = \m \left ( [0,\x] \right ) \;\; \left ( \text{respectively} \; \m_- (\x) = \m \left ( [-\x,0] \right )
\right )\,.
$$
To make $\m_{\pm}$  into a string, we need to define it over the whole of $[0,\infty]$.
It will suffice here to discuss  the extension of $\m_+$ into a string of length $\ell_+$;
the extension of $\m_-$ into a string of length $\ell_-$ is analogous and will be immediately obvious.

If $r$ is singular,
then it must be that $\r + \m \left ( [0,\r] \right )=\infty$. This forces $\ell_+  := \r$  and
$$
\m_+(\x) = \infty \;\;\text{for $\x \ge \ell_+$} \,.
$$
On the other hand, if $r$ is non-singular, we need to bring in the boundary condition satisfied by $\phi_{+} (\cdot,\lambda)$:
\begin{enumerate}
\item If $r$ belongs to $I$, we use the Feller-type boundary condition, i.e. we 
set $\ell_+ :=   \r +\frac{\beta_+}{\alpha_+}$ and
$$
\m_+(\x) =  \begin{cases}
\m \left ( [0,\r] \right ) + \frac{\gamma_+}{\beta_+} & \text{for $\r \le \x <\ell_+$} \\
\infty & \text{for $\x \ge \ell_+$}
\end{cases}\,.
$$

\item If $r$ does not belong to $I$, we
use the killing boundary condition, i.e. we set $\ell_+=\r$ and
$$
\m_+(\x) = \infty \;\; \text{if $\x \ge \ell_+$}\,.
$$
\end{enumerate}

Now, denote by $\left ( \xi_{\pm},\eta_{\pm} \right )$ the pair of fundamental solutions 
corresponding to the string $\m_{\pm}$ thus defined,
and by ${S}_{\pm}$ its characteristic function.
Set
$$
\xi (\x,\lambda) := \begin{cases}
\xi_-(-\x,\lambda) & \text{if $-\ell_- < \x \le 0$} \\
\xi_+(\x,\lambda) & \text{if $0 \le \x < \ell_+$}
\end{cases}
\;\;\text{and}\;\;
\eta (\x,\lambda) := \begin{cases}
-\eta_-(-\x,\lambda) & \text{if $-\ell_- < \x \le 0$} \\
\eta_+(\x,\lambda) & \text{if $0 \le \x < \ell_+$}
\end{cases}\,.
$$
Then the non-negative monotonic solution $\phi_{\pm}(\cdot,\lambda)$ of Equation (\ref{diffusionEquation}), suitably normalised, is obtained by setting 
$$
\phi_{\pm} (x,\lambda) = \xi (\x,\lambda) \mp \frac{1}{{S}_{\pm}(\lambda)} \eta (\x,\lambda)
$$
and we deduce easily that
\begin{equation}
\mp U_{\pm} (x_0,\lambda)  = \frac{s'(x_0)}{{S}_{\pm}(\lambda)} =  \lambda s'(x_0) {S}_{\pm}^\ast (\lambda)\,.
\label{stieltjesRiccatiRelation}
\end{equation}
It is clear from the foregoing construction
that $c_{\pm}^\ast = 0$; so this formula is equivalent to 
our earlier equation (\ref{knightFormula}), from which (\ref{stieltjesMeasure}) follows immediately.

\begin{remark}
In the case where the process $X$ is started at a reflecting boundary, then one of ${\m}_+$ or ${\m_-}$ is the zero string, and so
the corresponding generalised diffusion process is a one-sided diffusion.
\label{oneSidedRemark}
\end{remark}

\section{Proof of Theorem 1}
\label{proofOfTheorem1Section}

\begin{proof}
Suppose that $\sigma_{\pm}(x,\cdot)$ has infinitely many points of growth and that
the moment condition (\ref{momentCondition}) holds. Define a measure $\sigma$ on ${\mathbb R}_+$ by
\begin{equation}
\sigma (A) = \int_{\tilde{A}} z^{-1} \sigma_{\pm} (x,\d z)\,,\;\;\tilde{A} := \left \{ 1/z\,:\; z \in A \right \}\,,
\label{otherMeasure}
\end{equation}
Then $\sigma$ has infinitely many points of growth and 
satisfies the Stieltjes moment condition (\ref{stieltjesMomentCondition}). So we can write
$$
\int_0^{\infty} \frac{\sigma(\d z)}{\omega+z} = \frac{1}{\omega} \int_0^{\infty} \frac{\sigma(\d z)}{1+z/\omega} \sim
\sum_{j=0}^\infty \left [ \int_0^\infty (-z)^j \sigma (\d z) \right ] \omega^{-j-1}\;\;\text{as $\omega \rightarrow +\infty$}\,.
$$
Following Stieltjes \cite{St} (see also \cite{Ak,KW,NS}), one can construct from this series two sequences $\{ m_k\}$ and
$\{\ell_k\}$ of positive numbers such that, for every $n \in {\mathbb N}$,
$$
\cfrac{1}{m_0 \omega + \cfrac{1}{\ell_1+\cfrac{1}{m_1 \omega+\cfrac{1}{\ell_2+\cdots+\cfrac{1}{r_n}}}}}
= \sum_{j=0}^{n'} \left [ \int_0^\infty (-z)^j \sigma (\d z) \right ] \omega^{-j-1} + O(\omega^{-n'-2})
$$
as $\omega \rightarrow +\infty$, where
$$
n' = \begin{cases}
2n-1 & \text{if $r_n = \ell_n$} \\
2n & \text{if $r_n = m_n \omega$}
\end{cases}\,.
$$
Thus, in the same limit, we can write
\begin{multline}
\notag
\cfrac{1}{m_0 \omega + \cfrac{1}{\ell_1+\cfrac{1}{m_1 \omega+\cfrac{1}{\ell_2+\cdots+\cfrac{1}{r_n}}}}} \;+ O(\omega^{-n'-2})
= \int_0^\infty \frac{\sigma(\d z)}{\omega+z} \\
\overset{\underset{\downarrow}{z \rightarrow 1/z}}{=} \int_0^\infty \frac{\sigma_{\pm}(x,\d z)}{1+ z \omega}
= 1/\omega
\int_0^\infty \frac{\sigma_{\pm}(x,\d z)}{1/\omega+ z} = 1/2 \left [ \mp U_{\pm} (x,1/\omega) \pm U_{\pm}(x,0) \right ]\,,
\end{multline}
where we have made use of Equation (\ref{knightFormula}) to obtain the last equality. 
The first statement in the theorem follows if we take, for $n=0,\,1,\,\ldots$,
\begin{equation}
\mp u_{2n+1,\pm} (x) = \frac{m_n}{a(x)} \;\;\text{and}\;\; \mp u_{2 n+2,\pm} (x) = 2 \,\ell_n\,. 
\label{coefficientFormula}
\end{equation}

To prove the second statement, construct from the given $u_{n,\pm}(x)$
two sequences $\{m_n\}$ and $\{\ell_n\}$ of positive numbers via Equation
(\ref{coefficientFormula}). The hypothesis implies that the expansion
$$
\cfrac{1}{m_0 \omega + \cfrac{1}{\ell_1+\cfrac{1}{m_1 \omega + \cfrac{1}{\ell_2 + \cdots}}}}
$$
is a well-defined function of $\omega \in {\mathbb C} \backslash {\mathbb R}_-$. As shown by Stieltjes \cite{St},
since
$$
\sum_{n=0}^\infty m_n \quad \text{or} \quad \sum_{n=1}^{\infty} \ell_n
$$
diverges,
there is one and only one
measure $\sigma$ on ${\mathbb R}_+$ such that
$$
\int_0^\infty \frac{\sigma(\d z)}{\omega+z} = \cfrac{1}{m_0 \omega + \cfrac{1}{\ell_1+\cfrac{1}{m_1 \omega + \cfrac{1}{\ell_2 + \cdots}}}}\,.
$$
Furthermore, $\sigma$ satisfies the moment condition (\ref{stieltjesMomentCondition}). The uniqueness of $\sigma$
implies easily that Equation (\ref{otherMeasure}) holds, and so $\sigma_{\pm}(x,\cdot)$ satisfies the condition (\ref{momentCondition}).
The fact that it has infinitely many points of growth follows from the fact that the sequence of the $u_{n,\pm}(x)$ is infinite.
\end{proof}

\section{The expansion algorithm}
\label{expansionSection}

We now describe a method 
of obtaining the continued fraction which is based on the near-invariance of the Riccati equation under a certain linear fractional transformation.
The same idea, which can be traced back to Euler \cite{Eu}, has been used and developed more recently in a different context by Common \& Roberts \cite{CR}. Set
\begin{equation}
U(x,\lambda) =:  U_0(x,\lambda) =  u_0(x) + \frac{2 \lambda/a(x)}{U_1(x,\lambda)}\,,
\label{firstStep}
\end{equation}
where $u_0$ and $U_1$ are some functions which we shall specify presently. Substitution in Equation (\ref{riccatiEquation}) yields
\begin{multline}
\notag
u_0' - \frac{2 \lambda/a}{U_1} \frac{a'}{a} - \frac{2 \lambda/a}{U_1^2} U_1' + u_0^2 + 2 u_0 \frac{2 \lambda/a}{U_1} + 2 \lambda/a \frac{2 \lambda/a}{U_1^2} \\
+ 2 W' u_0 + 2 W' \frac{2 \lambda/a}{U_1} = 2 \lambda/a\,.
\end{multline}
Now {\em choose} $u_0$ so that it solves the {\em homogeneous} Riccati equation 
$$
u_0' + u_0^2 + 2 W' u_0 = 0\,.
$$
Then the equation satisfied by $U_1$ is
$$
U_1' + U_1^2 + 2 \left (  \frac{a'}{2 a} - u_0 - W' \right ) U_1 = 2 \lambda/a\,.
$$
This is of the same form as Equation (\ref{riccatiEquation}), save that $W'$ has been replaced by $a'/(2a)-u_0-W'$.
By iterating, we deduce the expansion
\begin{equation}
U (x,\lambda) = u_0(x) + \cfrac{2 \lambda/a(x)}{u_1(x) + \cdots+ \cfrac{2 \lambda/a(x)}{u_n(x)+\cfrac{2 \lambda/a(x)}{U_{n+1}(x,\lambda)}}}\,.
\label{continuedFractionExpansionWithRemainder}
\end{equation}
In this expression, $u_n$ satisfies the homogeneous Riccati equation
\begin{equation}
\frac{\d u_n}{\d x} + u_n^2  + 2 W_n'(x)\,  u_n = 0\,,
\label{coefficient}
\end{equation}
the remainder $U_n$ satisfies the {\em inhomogeneous} Riccati equation
\begin{equation}
\frac{\d U_n}{\d x} + U_n^2 + 2 W_n'(x)\, U_n = 2 \lambda/a
\label{remainder}
\end{equation}
and
$$
W_{n}' = \frac{a'}{2a} - u_{n-1} - W_{n-1}'\,,\;\; W_0' := W'\,.
$$
Theorem \ref{expansionTheorem} is thus proved.

The homogeneous Riccati equation always admits the trivial solution.
The calculation of the non-trivial solutions is sometimes
facilitated
by the following
\begin{proposition} 
The nontrivial solution $u_n$ of the homogeneous Riccati
equation of index $n$ satisfies the recurrence formula: for every $n \ge k \ge 0$,
$$
u_{n} (x) = \frac{\text{\em d}}{\text{\em d} x} \ln \left | \int \frac{\text{\em e}^{-2 W_k (x)}\,\text{\em d} x}{\prod_{j=k}^{n-1} \left ( a u_j^2 \right )}  \right |\,.
$$
\label{coefficientsProposition}
\end{proposition}

\begin{proof}
Since, by assumption, $u_j \ne 0$ for $k \le j \le n$, we deduce from the homogeneous equation that
$$
\frac{\d}{\d x} \frac{1}{u_j} - 2 W_j' \frac{1}{u_j} = 1\,.
$$
Hence, for $k \le j \le n$,
$$
u_j (x) = \frac{\e^{-2 W_j}}{\int \e^{-2 W_j}\,\d x} = \frac{\d}{\d x} \ln \left | \int \e^{-2 W_j} \,\d x \right |\,.
$$

Let us write
$$
W_n (x) = \int W_n'(x) \,\text{\em d} x\,,
$$
where the indefinite integral means that $W_n$ is any primitive of the integrand. 
For $n > k$,
\begin{multline}
\notag
-2 W_n = -2 \int W_n'(x)\,\d x = -2 \int \left [ \frac{1}{2} \frac{a'}{a} - u_{n-1} - W_{n-1}' \right ]\,\d x \\
= -\ln a + 2 W_{n-1} + 2 \ln  \left | \int \e^{-2 W_{n-1}} \,\d x \right |
\end{multline}
and so
$$
\e^{-2 W_n} = \frac{1}{a} \e^{2 W_{n-1}} \left | \int \e^{-2 W_{n-1}} \,\d x  \right |^2 = \frac{\e^{-2 W_{n-1}}}{a u_{n-1}^2}\,.
$$
The required formula follows easily by iterating.
\end{proof}

\section{A probabilistic interpretation of the algorithm}
\label{probabilisticSection}
It is straightforward to give a probabilistic interpretation
of the first few coefficients in the continued fraction. Indeed,
from Equations (\ref{levyKhintchineFormula}-\ref{pitmanYorFormula}) and
Equation (15) of \cite{PY},
$$
u_{0,-}(x) - u_{0,+}(x)
$$  
is inversely proportional to the mean of the local time spent at $x$. Also,
$$
\frac{1}{\mp a(x) u_{1,\pm}(x)} = \int_0^\infty y \nu_{\pm} (x,\d y)
$$
and hence $\mp u_{1,\pm}(x)$
is inversely proportional to the average duration of the finite excursions to the right ($+$) or left ($-$).

Our aim in this section is to gain some insight into the probabilistic content of
the {\em expansion  algorithm} itself. If the initial diffusion $X$ has a Riccati variable $U$ with a continued fraction expansion, then, as can be seen from
the proof of Theorem \ref{expansionTheorem} given in the previous section, the algorithm
produces a sequence of
$W_n$ and a sequence of remainders $U_n$. Each $U_n$ can be thought of as a Riccati variable of some diffusion, say $X_n$, corresponding to the ``environment'' $W_n$,
and it is then
natural to investigate the relationship between adjacent diffusions in this sequence.  
In pursuing this line of thought, it is important to bear in mind the following points:
\begin{enumerate}
\item The $W_n$ depend on which of the two Riccati variables $U_+$ and $U_-$ is being expanded.
In particular, if {\em both} $U_+$ and $U_-$ have a continued fraction expansion, then there are {\em two} sequences of $W_n$.
\item It takes {\em two} strings to specify a regular diffusion--- each string corresponding to a one-sided diffusion (see Definition
\ref{oneSidedDefinition} and Remark \ref{oneSidedRemark}). Each Riccati variable specifies
one string, and the regular diffusion is obtained by welding the strings
together.
\end{enumerate}
In the remainder of this section, we will show how, from {\em one} continued fraction expansion of a Riccati variable of the regular diffusion $X$, one can
construct a particular sequence of diffusions $X_n$ such that the remainder $U_n$ is a Riccati variable of $X_n$. 
It should be clear from the points
just made that sequences with this property cannot be unique. 
Our particular construction is inspired by the 
following key formula, proved by
Pitman \& Yor \cite{PY}:
\begin{equation}
\psi_{\pm}(x,\lambda) 
= \psi_{\pm} (x,0) \mp \frac{1}{2} \frac{\d}{\d y} \Bigl |_{y=x \pm} {\mathbb E}_y \left ( \e^{- \lambda H (x) } \Bigl |  \,H (x) < \infty \right )\,.
\label{keyPitmanYorFormula}
\end{equation}
To simplify matters, we assume that 
the speed measure $m$ of $X$ is absolutely
continuous with respect to the Lebesgue measure and, with some abuse of notation, write
$$
m(A) = \int_A m(x)\,\d x
$$
for every measurable set $A$.
The instantaneous variance $a$ and the instantaneous drift $b$ of $X$ are then given
by
$$
a = \frac{2}{m s'} \;\;\text{and}\;\; b = \frac{1}{m} \frac{\d}{\d x} \frac{1}{s'}\,.
$$

\subsection{The $h$-transform}
\label{hTransformSubsection}

For definiteness, let us suppose that we work with a continued fraction expansion of the Riccati variable $U_+$. Taking the $+$ sign in Equation
(\ref{keyPitmanYorFormula}) and making use of the identity (\ref{pitmanYorFormula}), we obtain
$$
U_{+}(x,\lambda) 
= U_{+} (x,0) - \frac{\d}{\d y} \Bigl |_{y=x+} {\mathbb E}_y \left ( \e^{- \lambda H (x) } \Bigl |  \,H (x) < \infty \right )\,.
$$
We will show that the second term on the right-hand side is the Riccati variable of a diffusion, say $Y$,  obtained from
$X$ by conditioning, and we will express the characteristics of $Y$ in terms of the characteristics of $X$. In particular,
it will follow from our construction that
the Riccati variables $U_{\pm}^Y$ of $Y$ are related to those of $X$ via
\begin{equation}
U_{+} (x,\lambda) = U_{+} (x,0) + U_{+}^Y(x,\lambda) \;\; \text{and}\;\; U_{-} (x,\lambda) = U_{+} (x,0) + U_{-}^Y(x,\lambda)  \,.
\label{riccatiVariableOfY}
\end{equation}

If $U_{+}(\cdot,0)$ vanishes, then we take $Y=X$. Otherwise,
$X$ is transient and the event
$$
\lim_{t \rightarrow \zeta} X_t = r\,,
$$
where $\zeta$ denotes the lifetime of $X$, occurs with positive probability. Let $l < z < x$ and denote by $Y$ the process obtained from $X$ by conditioning on the event
$$
H(z) < \infty\,.
$$ 
We have (see \cite{BS}, II.12)
$$
{\mathbb P}_x \left ( H(z) < \infty \right ) = \phi_+(x,0) =: h(x)\,,
$$
where ${\mathbb P}_x$ is the probability measure associated with $X$ started at $x$ and $\phi_+$ is normalised
so that $\phi_+(z,0)=1$. 

For $z < y < x$, denote by $d y$ an interval of infinitesimal length centered on $y$. Then
\begin{multline}
\notag
{\mathbb P}_x \left ( Y_t \in d y \right ) =  {\mathbb P}_x \left ( X_t \in d y \Bigl | \,H(z) < \infty \right )
= \frac{{\mathbb P}_x \left (  X_t \in d y,\,H(z) < \infty  \right )}{{\mathbb P}_x \left ( H(z) < \infty \right )} \\
= \frac{{\mathbb P}_x \left (  X_t \in d y \right ) {\mathbb P}_y \left (  H(z) < \infty \right )}{{\mathbb P}_x \left ( H(z) < \infty \right )}
= \frac{{\mathbb P}_x \left (  X_t \in d y \right ) h(y)}{h(x)}\,.
\end{multline}
Since $h$ is excessive (see \cite{BS}, II.30), this calculation shows that $Y$ is the $h$-transform of $X$.
This result does not depend on the
particular choice of $z$ as long as $l < z < x$. 
Clearly,
$$
\phi_{+}^Y = \frac{\phi_{+}}{h}\,.
$$
The characteristics of $Y$
are displayed in Table \ref{transformationTable}.

The foregoing discussion assumed that one is working with the Riccati variable $U_+$. If, instead, one is working with the
other Riccati variable $U_-$, then one should use $h = \phi_-(\cdot,\lambda)$ to define $Y$. It is then easily verified 
that
$$
U_{+} (x,\lambda) = U_{-} (x,0) + U_{+}^Y(x,\lambda) \;\; \text{and}\;\; U_{-} (x,\lambda) = U_{-} (x,0) + U_{-}^Y(x,\lambda) 
$$
holds instead of Equation (\ref{riccatiVariableOfY}).

\begin{table}
\begin{tabular}{c c c c c}
$X$ & $\xrightarrow[]{\text{$h$-transform}}$ & $Y$ & $\xrightarrow[]{\text{Krein duality}}$ & ${\mathscr T}_h (X)$ \\
\hline
\hline
$m$ & & $h^2 m$ & & $h^{-2} s'$  \\
$s'$ & & $h^{-2} s'$ & & $h^2 m$ \\
$a$ & & $a$ & & $a$ \\
$b$ & & $b + a h'/h$ & & $a'/2 - b - a h'/h$ \\
$U_{\pm}(\cdot,\lambda)$ & & $U_{\pm}(\cdot,\lambda)-h'/h$ & & $\frac{2 \lambda/a}{U_{\pm}(\cdot,\lambda)-h'/h}$
\end{tabular}
\caption{The transformation ${\mathscr T}_h$ and its effect on the characteristics of the diffusion.}
\label{transformationTable}
\end{table}

\subsection{Krein duality}
\label{dualitySubsection}
Given Equation (\ref{riccatiVariableOfY}), we ask next
for a diffusion, say $Z$, whose Riccati variables are 
related to those of $Y$ via
\begin{equation}
U_{\pm}^Z (x,\lambda)  \,U_{\pm}^Y (x,\lambda)=  2 \lambda/a(x)\,.
\label{riccatiVariableOfZ}
\end{equation} 
We emphasise that this equation should be understood as specifying {\em both} $U_+^Z$ and $U_-^Z$.
Equation (\ref{stieltjesRiccatiRelation}) suggests that the Krein dual of $Y$ might be a suitable candidate. This raises
some technical issues:
$Y$ is a well-defined diffusion process, but its behaviour at the boundary points may differ from that of $X$, and it is not
immediately clear whether $Y$ is non-singular or whether it is even a generalised diffusion process in the sense
of \S \ref{gapSubsection}. To avoid these complications, we shall
be content to {\em suppose} that $Y$ belongs to the particular class of non-singular diffusions described in the statement of the following
proposition, so that its Krein dual exists and is a non-singular diffusion in the same class.
The resulting concept of duality is then analogous to that used by Jansons \cite{Ja}, Soucaliuc \cite{So} and T\'{o}th \cite{To} (see the definition
of {\em conjugate diffusion} in his Appendix 1).

\begin{proposition}
Let $X$ be a non-singular diffusion process, started at $x$, on an interval $I$ that includes no non-singular boundary points--- except possibly reflecting ones.
Let $X^{\ast}$ be the Krein dual of $X$, in the same scale as $X$.
Then
\begin{enumerate}
\item $X^{\ast}$ is a non-singular diffusion process, started at $x$, on an interval $I^{\ast}$ with the same endpoints as $I$
and that includes no non-singular boundary points--- except
possibly reflecting ones;
\item the speed measure of $X$ is the scale function of $X^{\ast}$ and vice-versa;
\item the Riccati variables $U_{+}^{\ast}$ and $U_{-}^{\ast}$ of $X^{\ast}$ are related to the Riccati variables $U_{+}$ 
and $U_{-}$ of $X$ via
$$
U_{\pm} (x,\lambda) U_{\pm}^{\ast} (x,\lambda) =  2 \lambda/a(x)\,,
$$
where $a$ is the common infinitesimal variance of $X$ and $X^\ast$.
\end{enumerate}
\label{dualProposition}
\end{proposition}

\begin{proof}
We work ``backwards'', i.e. from $X$, we define a non-singular process $\hat{X}$ on an interval $I^\ast$ with the same endpoints as $I$, started
at $x$.
Then we show that its
associated pair of strings is the same as that of the generalised diffusion process $X^{\ast}$. For the sake of convenience, we shall
use in the proof the same notation
for the process $X$ as in \S \ref{diffusionSubsection}. 

In the interior of $I^\ast$, the process $\hat{X}$ is determined by the generator
$$
{\mathscr G}_W^{\ast} := \frac{1}{s'(x)} \frac{\d}{\d x} \left ( \frac{1}{m(x)} \frac{\d}{\d x} \right )\,.
$$ 
The boundary behaviour of $\hat{X}$ is as follows: for every non-singular $p \in \{ l,r \}$, if $p$ is reflecting for $X$ then it is killing for
$\hat{X}$ and $p \notin I^\ast$; if $p$ is killing for $X$, then $p$ is reflecting for $\hat{X}$ and $p \in I^\ast$.

We remark that, for $\l < \x <\r$,
$$
{\mathscr G}_W^{\ast} \overset{\underset{\downarrow}{\x=s(x)}}{=} \frac{\d}{\d \x} \left ( \frac{s'(x)}{m(x)} \frac{\d}{\d \x} \right ) = \frac{\d^2}{\d \x \d \m (\x)}\,,
$$
where $\d \m (\x) = \m(\d \x) = m(\d x)$. Recalling the definition of the string $\m_{\pm}$, this shows that the natural scale of $\hat{X}$ is
$$
\y = \m(\x) = \begin{cases}
\m_+ (\x) & \text{if $0 < \x < \r$} \\
-\m_-(-\x) & \text{if $\l < \x \le 0$}
\end{cases}
$$
and we deduce
\begin{equation}
\hat{\m}_+ (\y) = \m_+^\ast (\y)\,, \; 0 < \y < \m_+(\r)\,, \quad \hat{\m}_- (\y) = \m_-^\ast (\y)\,, \; 0 \le \y < \m_-(-\l)\,.
\label{twoEqualities}
\end{equation}

To show that $\hat{X} = X^\ast$, there only remains to extend the validity of these two equalities to every $\y \in [0,\infty]$. 
Consider the extension of the first equality.
As explained in \S \ref{diffusionSubsection}, the definition of $\hat{\m}_+$ in $[\m_+(\r),\infty]$ depends on the nature of the boundary point $r$.
We bear in mind that,
since the speed measure of $\hat{X}$ is the scale function of $X$ and vice-versa, $r$ is entrance for $\hat{X}$ if and only if it is exit for $X$ and
vice-versa.
There are three cases to consider: 
\begin{enumerate}[(i)]
\item If $r$ is singular for $\hat{X}$, then $\hat{\m}_+(\y) = \infty$ for $\y \ge \m_+(\r)$. Since $r$ is also singular for $X$, we have
$\r + \m_+(\r) = \infty$. If $\m_+(\r) = \infty$, there is nothing to extend; if $\m_+(\r) < \infty$, then $\r = \infty$ and, 
since $\m_+^\ast$ is by definition the right-continuous inverse of $\m_+$,
$$
\m_+^\ast (\y) = \infty = \hat{\m}_+(\y) \;\;\text{for $\y \ge \m_+(\r)$}\,. 
$$

\item If $r$ is reflecting for $\hat{X}$, then
$$
\hat{\m}_+(\y) = \begin{cases}
\r & \text{for $\m_+(\r) \le \y < \infty$} \\
\infty & \text{for $\y = \infty$}
\end{cases}\,.
$$
On the other hand, $r$ is killing for $X$, i.e. $\m_+(\x) = \infty$ for $\x \ge \r$. Again, since $\m_+^\ast$ is by definition the right-continuous
inverse of $\m_+$, we deduce $\m_+^\ast (\y) = \hat{\m}_+(\y)$ for $y \ge \m_+(\r)$.

\item If $r$ is killing for $\hat{X}$, then it is reflecting for $X$, and the desired result follows from (ii) by symmetry.
\end{enumerate}
The extension of the second equality in (\ref{twoEqualities}) follows along the same lines. 
We have thus shown that $\hat{X}=X^{\ast}$. 

The first two statements in the proposition follow immediately.
For the last statement, we apply Equation (\ref{stieltjesRiccatiRelation}) to $X$ and to $X^{\ast}$; this gives
$$
U_{\pm} (x,\lambda)\,U_{\pm}^{\ast} (x,\lambda) = \frac{s'(x)}{S_{\pm}(\lambda)} \frac{m(x)}{S_{\pm}^\ast (\lambda)} 
= \frac{2/a(x)}{1/\lambda}\,.
$$
\end{proof} 
Table \ref{transformationTable} shows how the characteristics of a process transform under Krein duality 
in the situation envisaged by Proposition \ref{dualProposition}.

To summarise the foregoing discussion, one can--- at least in some cases--- describe the first iteration of the expansion algorithm in terms
of a map ${\mathscr T}_h$ obtained by composing two transformations:
$$
X \xrightarrow[]{\text{$h$-transform}} Y \xrightarrow[]{\text{Krein duality}} Z =: {\mathscr T}_h \left (X \right )\,.
$$
More precisely, expanding the Riccati variable $U_{+}$ of $X$ leads, after one step, to the Riccati variable of the ``new'' diffusion ${\mathscr T}_{h}(X)$ 
where $h = \phi_{+}(\cdot,0)$. 

If, instead, one works with the Riccati variable $U_-$, then 
$h = \phi_{-}(\cdot,0)$ is used to construct the next diffusion 
${\mathscr T}_h \left (X \right )$.

\subsection{The Ciesielski--Taylor theorem}
\label{bianeCiesielskiTaylorSubsection}
Ciesielski and Taylor \cite{CT} noticed that the total time spent by a $(d+2)$-dimensional standard Brownian motion inside the unit ball in ${\mathbb R}^{d+2}$
has the
same distribution as the first hitting time of the unit ball in ${\mathbb R}^{d}$ by a $d$-dimensional standard Brownian motion. This result can be expressed
in terms of one-dimensional diffusions as an identity in law between some occupation time of a Bessel process (of parameter $p$) and a hitting time of another Bessel process (of parameter $p-1$). 
Biane \cite{Bi} generalised the Ciesielski--Taylor theorem to other pairs of diffusions. As we shall see presently, 
his construction of these pairs can be described
in terms of the map ${\mathscr T}_h$ defined in the previous section. A similar observation has already been made by T\'{o}th in \cite{To}.

Suppose that
(1) $s(l) = -\infty$; (2) $s(r)=0$; (3) if $r<\infty$, $r$ is killing, and
$$
\text{(4)} \;\;
\begin{cases}
s(y) \,m(y) \xrightarrow[y \rightarrow l]{} 0 & \text{if $l > -\infty$} \\
 & \\
\int_{-\infty}^x s^2 (y) \,m (\text{\em d} y) = \infty & \text{if $l=-\infty$}
\end{cases}\,.
$$
The hypothesis ensures that
$$
X_t \xrightarrow[t \rightarrow \zeta]{} r \;\;\text{almost surely}\,.
$$
Then $h := \phi_+(\cdot,0)= -s > 0$, and the $h$-transform $Y$ of $X$
is a non-singular diffusion such that $r$ is singular and, if $l$ is non-singular, then it is killing. Its dual
$Y^\ast = \mathscr{T}_h(X)$ is therefore a well-defined diffusion.

We seek a diffusion $Z$, with a suitable initial law, such that, for every $l < x \le y < r$,
\begin{equation}
\inf \left \{ t \ge 0 \,:\; Z_t = y \right \} \overset{\text{(law)}}{=}  \int_0^{\zeta} {\mathbf 1}_{\{ X_{t} \le y \}} \,\d t \,.
\label{ciesielskiTaylor}
\end{equation}
Following Biane \cite{Bi}, the Feynman--Kac formula gives
\begin{equation}
{\mathbb E}_x \left ( \exp \left [ -\lambda \int_0^{\zeta} {\mathbf 1}_{\{ X_{t} \le y \}} \,\d t \right ] \right )
= \frac{\phi_-(x,\lambda)}{\phi_-(y,\lambda)-\frac{h(y)}{h'(y)} \phi_-'(y,\lambda)}\,.
\label{feynmanKac}
\end{equation}
The task is therefore to find $Z$ such that
\begin{equation}
\phi_-^Z (y,\lambda) = \phi_-(y,\lambda)-\frac{h(y)}{h'(y)} \phi_-'(y,\lambda)
\label{bianeEquation1}
\end{equation}
and a distribution $\nu$ for $Z_0$ such that
\begin{equation}
\int_{l}^x \phi_-^Z(z,\lambda) \,\nu (\d z) = \phi_-(x,\lambda)\,.
\label{bianeEquation2}
\end{equation}
To find the former, we start with Biane's observation that Equation (\ref{bianeEquation1}) implies
$$
\frac{\d}{\d y} \phi_-^Z (y,\lambda) = -s(y) \frac{\d}{\d y} \left [ \frac{1}{s'(y)} \frac{\d}{\d y} \phi_- (y,\lambda) \right ]\,.
$$
So the ``Riccati version'' of Equation (\ref{bianeEquation1}) is
\begin{multline}
\notag
U_{-}^Z (y,\lambda) = \frac{-s(y) \frac{\d}{\d y} \left [ \frac{1}{s'(y)} \frac{\d}{\d y} \phi_- (y,\lambda) \right ]}{\phi_-(y,\lambda)-\frac{h(y)}{h'(y)} \phi_-'(y,\lambda)} 
= -m(y) s(y) \frac{{\mathscr G}_W \phi_-(y,\lambda)}{\phi_-(y,\lambda)-\frac{h(y)}{h'(y)} \phi_-'(y,\lambda)} \\
= -\frac{s(y)}{s'(y)} \frac{2/a(y)\, \lambda\,\phi_-(y,\lambda)}{\phi_-(y,\lambda)-\frac{h(y)}{h'(y)} \phi_-'(y,\lambda)}
= \frac{2 \lambda/a(y)}{U_-(y,\lambda) - \frac{h'(y)}{h(y)}}\,.
\end{multline}
Table \ref{transformationTable} then shows that, by taking $Z = {\mathscr T}_h \left (X \right )$,
we obtain a diffusion such that this does indeed hold, and hence also Equation (\ref{ciesielskiTaylor})--- given the right $\nu$.

\section{Some deterministic examples}
\label{exampleSection}
There are two methods for calculating the continued fraction
coefficients: the first uses the algorithm of \S \ref{expansionSection}; the second proceeds
by tracking the diffusions produced by the recurrence
\begin{equation*}
X \xrightarrow[]{}  {\mathscr T}_{h_0} \left ( X \right) \xrightarrow[]{} {\mathscr T}_{h_1} \circ {\mathscr T}_{h_0} \left ( X \right) \xrightarrow[]{} \hdots
\,.
\end{equation*}
The latter method provides greater insight but is not always applicable; the former method can always be used, but the determination
of the constants of integration can be tedious.
In this section, we study some well-known diffusions in deterministic environments which will serve to illustrate the two approaches.
The processes that we consider have the remarkable
property of belonging to a parametrised class that is closed under 
both the $\phi_{\pm}(\cdot,0)$-transform and Krein duality; see Table \ref{diffusionTable}. 
Our main reference is \cite{BS}, Appendix I. 

\subsection{Brownian motion with drift}
\label{brownianMotionWithDriftSubsection}
Let $B^{\mu}$ be Brownian motion with drift $\mu$ and set $X = B^{\mu}$. The generator is
$$
{\mathscr G}_W = \frac{1}{2} \frac{\d^2}{\d x^2} + \mu \frac{\d}{\d x}
$$
and the state space is ${\mathbb R}$; the endpoints are natural and so there are no boundary conditions. We deduce that the Riccati equation has a unique
positive solution and a unique negative solution. We find
$$
\phi_{\pm} = \exp \left [ x \left ( -\mu \mp \sqrt{\mu^2+2 \lambda} \right ) \right ]\;\;\text{and so}\;\;
U_{\pm} = -\mu \mp \sqrt{\mu^2+2 \lambda}\,.
$$
The Stieltjes--Perron inversion formula yields
$$
\sigma_{\pm} (x,\,\d z) = \frac{1}{\pi} \frac{\sqrt{2 z-\mu^2}}{2 z}\,\d z\,, \;\; z > \mu^2/2
$$
and
$$
\sigma_{\pm} (x,\{0\}) = \frac{1}{2} \left [ |\mu| \pm \mu \right ]\,.
$$
Condition (\ref{momentCondition}) is therefore satisfied unless $\mu=0$ (standard Brownian motion). 
In fact, the expansion of $U_{\pm}$ is easily obtained by elementary means: assuming for definiteness that $\mu > 0$, we have
\begin{multline}
\notag
U_- = \sqrt{\mu^2+2\lambda}-\mu = \frac{2 \lambda}{\sqrt{\mu^2+2\lambda}+\mu} \\ = \frac{2 \lambda}{2 \mu + U_-} 
= \cfrac{2 \lambda}{2 \mu + \cfrac{2 \lambda}{2 \mu + \cfrac{2 \lambda}{2 \mu + \cdots}}}
\end{multline}
and
$$
-U_+ = \sqrt{\mu^2+2\lambda}+\mu = 2\mu + U_- = 2 \mu + \cfrac{2 \lambda}{2 \mu + \cfrac{2 \lambda}{2 \mu + \cfrac{2 \lambda}{2 \mu + \cdots}}}\,.
$$

Let us demonstrate how these expansions can be obtained by the algorithm of \S \ref{expansionSection}.
Again, for definiteness, suppose that $\mu > 0$. We need to distinguish two cases, namely $u_0=0$ and $u_0 \ne 0$. 

Take $u_0 = 0$. Then $b_1 = -\mu$, $W_1 = -\mu x$ and so
$$
u_1 (x) = \frac{\e^{2 \mu x}}{\frac{\e^{2 \mu x}}{2 \mu}+c}\,,
$$
where $c$ is a constant of integration. There is no possibility of choosing this constant so that $u_1$ is negative. But by taking $c \ge 0$,
we obtain $u_1 > 0$.
Furthermore,
$$
u_2(x) = \frac{\e^{2 \mu x}/u_1^2}{\int \e^{2 \mu x}/u_1^2\,\d x} = 
\frac{\frac{1}{4 \mu^2} \e^{2 \mu x} + \frac{c}{\mu}+ c^2 \e^{- 2\mu x}}{\frac{1}{8 \mu^3} \e^{2 \mu x} + \frac{c}{\mu} x - \frac{c^2}{2 \mu} \e^{-2 \mu x}+C}
$$
where $C$ is another constant of integration. The only way to ensure that $u_2$ is also positive is by taking $c=0$. Hence
$$
u_1 (x) = 2 \mu \;\;\text{and}\;\; u_2(x) = \frac{\frac{1}{4 \mu^2} \e^{2 \mu x}}{\frac{1}{8 \mu^3} \e^{2 \mu x}+C}\,.
$$
By iterating, we obtain the expansion of $U_-$.

The other possibility is to take $u_0 \ne 0$. Since $W_0 = \mu x$, we find
$$
u_0(x) = \frac{\e^{-2 \mu x}}{\frac{-1}{2 \mu} \e^{-2 \mu x}+c}\,,
$$
where $c$ is a constant of integration. No choice of $c$ can make $u_0$ positive but, by taking $c \le 0$, we can make it negative.
Then
$$
u_1 (x) = \frac{\e^{-2 \mu x}/u_0^2}{\int \e^{-2 \mu x}/u_0^2\,\d x}\,. 
$$
The numerator is positive; the denominator is
$$
- \frac{1}{8 \mu^3} \e^{-2 \mu x} - \frac{c}{\mu} x + \frac{c^2}{2 \mu} \e^{2 \mu x} + C\,,
$$
where $C$ is another constant of integration. This expression cannot be negative for every $x$ unless $c=0$. We deduce
$$
u_0(x) = -2 \mu \;\;\text{and}\;\;u_1(x) = \frac{\frac{1}{4 \mu^2} \e^{-2 \mu x}}{\frac{-1}{8 \mu^3} \e^{-2 \mu x}+C}\,.
$$
By iterating, we obtain the expansion of $U_+$.

Next, we demonstrate how the expansion may be interpreted in terms of a sequence of diffusions obtained via $h$-transforms and Krein duality. We have
$$
\phi_{\pm} (x,0) = \exp \left [ \left (-\mu \mp | \mu | \right ) x \right ]\,.
$$
Suppose that $\mu > 0$. Then $\phi_{-}(\cdot,0) = 1$, $\phi_+(x,0) = \e^{-2 \mu x}$ and, from Table \ref{transformationTable}, it is easily deduced that the $\phi_+(\cdot,0)$-transform
of $B^{\mu}$ is $B^{-\mu}$. On the other hand, if $\mu<0$, then $\phi_{+}(\cdot,0) = 1$, $\phi_-(x,0) = \e^{-2 \mu x}$, and so it is the $\phi_-(\cdot,0)$-transform of
$B^{\mu}$ that yields $B^{-\mu}$.
Also, from Table \ref{transformationTable}, it is immediate that the Krein dual of $B^{\mu}$ is $B^{-\mu}$. Putting these results together, we find, for example,
that the expansion of $U_+$ for $\mu >0$ corresponds to the sequence
\begin{multline}
\notag
X = B^{\mu} \xrightarrow[]{\text{$\phi_+(\cdot,0)$-transform}} B^{-\mu}  \xrightarrow[]{\text{Krein duality}} {\mathscr T}_{h_0} \left ( X \right ) =  B^{\mu}  \\
\xrightarrow[]{\text{$\phi_+(\cdot,0)$-transform}} B^{-\mu}  \xrightarrow[]{\text{Krein duality}}  {\mathscr T}_{h_1} \circ {\mathscr T}_{h_0} \left ( X \right ) = B^{\mu} \;\;\text{etc.} 
\end{multline}

\begin{table}
\begin{tabular}{l | l | l | l | l}
$X$ & Parameters & $X^-$ & $X^+$ & $X^\ast$ \\
\hline
\hline
$B^{\mu}$ & \text{\small $\mu < 0$} & $B^{-\mu}$ & $B^{\mu}$ &  $B^{-\mu}$ \\
\hline
$B^{\mu}$ & \text{\small $\mu > 0$} & $B^{\mu}$ & $B^{-\mu}$ &  $B^{-\mu}$ \\
\hline
$\text{BES}(p)$ & \text{\small $p \le -1$} & $\text{BES}(-p)$ & $\text{BES}(p)$ & $\text{BES}(-p-1)$ \\
\hline
$\text{BES}(p)$ & \text{\small $-1<p<0$,} & $\text{BES}(-p)$ & $\text{BES}(p)$, & $\text{BES}(-p-1)$, \\
 & \text{\small $0$ killing} & & \text{\small $0$ killing} &   \text{\small $0$ reflecting} \\
\hline
$\text{BES}(p)$ & \text{\small $-1<p<0$}, & $\text{BES}(p)$, & $\text{BES}(p)$, & $\text{BES}(-p-1)$,\\
 & \text{\small $0$ reflecting} & \text{\small $0$ reflecting} & \text{\small $0$ reflecting} & \text{\small $0$ killing} \\
\hline
$\text{BES}(p)$ & \text{\small $p=0$} & $\text{BES}(0)$ & $\text{BES}(0)$ & $\text{BES}(-1)$ \\
\hline
$\text{BES}(p)$ & \text{\small $0<p<1$} & $\text{BES}(p)$ & $\text{BES}(-p)$, & $\text{BES}(-p-1)$ \\
 & & &  \text{\small $0$ killing} & \\
\hline
$\text{BES}(p)$ & \text{\small $p \ge 1$} & $\text{BES}(p)$ & $\text{BES}(-p)$ & $\text{BES}(-p-1)$
\end{tabular}
\caption{The $\phi_{\pm}(\cdot,0)$-transform $X^{\pm}$ and the Krein dual $X^\ast$ for some processes $X$ in a deterministic environment: $B^{\mu}$ denotes Brownian motion
with drift $\mu$ and $\text{BES}(p)$ denotes a Bessel process of parameter $p$.}
\label{diffusionTable}
\end{table}

\subsection{A Bessel process}
\label{besselSubsection}
Let $\text{BES}(p)$ denote the Bessel process with parameter $p$ and, for $p > 0$, let $X=\text{BES}(p)$. Then
$$
{\mathscr G}_W = \frac{1}{2} \frac{\d^2}{\d x^2} + \frac{p+\frac{1}{2}}{x} \frac{\d}{\d x}\,, \;\; x>0\,,
$$
$0$ is entrance-not-exit, $\infty$ is natural,
$$
\phi_-(x,\lambda) = x^{-p} I_{p} \left ( \sqrt{2 \lambda} x \right )\,,\;\;\phi_+(x,\lambda) = x^{-p} K_{p} \left ( \sqrt{2 \lambda} x  \right )
$$
and
$$
U_- (x,\lambda) = \sqrt{2 \lambda} \frac{I_{p+1} (\sqrt{2 \lambda} x)}{I_p(\sqrt{2 \lambda} x)}\,,\;\; -U_+ (x,\lambda) = \sqrt{2 \lambda} \frac{K_{p+1} (\sqrt{2 \lambda} x)}{K_p(\sqrt{2 \lambda} x)}\,.
$$
To obtain these expressions for the Riccati solutions, we have made use of the three-term recurrence relations satisfied by the Bessel
functions $I_p$ and $K_p$:
$$
I_{p-1} (z) - I_{p+1}(z) = \frac{2 p}{z} I_p(z)\,,\;\; 2 I_p'(z) = I_{p-1}(z)+I_{p+1}(z)
$$
and
$$
K_{p-1} (z) - K_{p+1}(z) = -\frac{2 p}{z} K_p(z)\,,\;\; -2 K_p'(z) = K_{p-1}(z)+K_{p+1}(z)\,.
$$

The Stieltjes measures are \cite{IK}:
$$
\sigma_-(x,\{0\}) = 0\,, \quad \sigma_-(x,\d z) = \sum_{k=1}^\infty \frac{1}{x} \,\delta_{\frac{j_{p,k}^2}{2 x^2}}\,,
$$
where the $j_{p,k}$ are the positive zeroes of the Bessel function $J_p$, 
and
$$
\sigma_+(x,\{0\}) = \frac{p}{x}\,, \quad \sigma_{+}(x,\d z) = \frac{1}{\pi^2 z x} \frac{\d z}{J_p^2 (\sqrt{2 z}x) + Y_p^2 (\sqrt{2 z}x)}\,.
$$
Since the $j_{p,k}$ increase linearly with $k$ as $k \rightarrow \infty$, $\sigma_-(x)$ satisfies the moment condition. On the other hand, since
$$
\sigma_{+}(x,\d z) \sim c(p,x) \,z^{p-1} \,\d z \quad \text{as $z \rightarrow 0$}\,,
$$
the moment condition is not satisfied by $\sigma_+(x,\cdot)$. By using the second of the recurrence relations satisfied by the Bessel functions, it
is straightforward to verify that
\begin{multline}
\notag
U_- = \sqrt{2 \lambda} \,\frac{I_{p+1}(\sqrt{2 \lambda}x)}{I_{p}(\sqrt{2 \lambda} x)} = \frac{\sqrt{2 \lambda}}{\frac{2(p+1)}{\sqrt{2 \lambda} x}+\frac{I_{p+2}(\sqrt{2 \lambda} x)}{I_{p+1}(\sqrt{2 \lambda} x)}} \\ = \frac{2 \lambda}{\frac{2(p+1)}{x}+ \sqrt{2 \lambda} \,\frac{I_{p+2}(\sqrt{2 \lambda}  x)}{I_{p+1}(\sqrt{2 \lambda} x)}} 
= \cfrac{2 \lambda}{\frac{2(p+1)}{x}+\cfrac{2 \lambda}{\frac{2(p+2)}{x}+\cfrac{2 \lambda}{\frac{2(p+3)}{x}+\cdots}}}\,.
\end{multline}
Also, despite the fact that $\sigma_+(x)$ does not satisfy the moment condition,
\begin{multline}
\notag
-U_+ = \sqrt{2 \lambda} \frac{K_{p+1} \left ( \sqrt{2 \lambda} x \right )}{K_p \left ( \sqrt{2 \lambda} x \right )} = \frac{2p}{x} +
\frac{2 \lambda}{\sqrt{2 \lambda} \frac{K_{p} \left ( \sqrt{2 \lambda} x \right )}{K_{p-1} \left ( \sqrt{2 \lambda} x \right )}} \\
= \frac{2 p}{x} + \cfrac{2 \lambda}{\frac{2(p-1)}{x}+\cfrac{2 \lambda}{\frac{2(p-2)}{x}+\cfrac{2 \lambda}{\frac{2(p-3)}{x}+\cdots}}}\,.
\end{multline}
This expansion does not contradict Theorem \ref{stieltjesTheorem}, for the coefficients cannot all be positive.

It is also possible to obtain both expansions from the algorithm of \S \ref{expansionSection}. First, consider the choice $u_0 = 0$. Then $b_1 = -b$ and
$$
\e^{-2 W_1 (x)} = x^{2p+1}\,.
$$
We obtain
$$
u_1 (x) = \frac{x^{2p+1}}{\frac{x^{2p+2}}{2p+2}+c}
$$
for some constant of integration $c$. Then
$$
u_2 (x) = \frac{\e^{-2 W_1}/u_1^2}{\int \e^{-2 W_1}/u_1^2 \d x} = \frac{x^{2p+1} \left [ \frac{x}{2p+2} + c x^{-2p-1} \right ]^2}{\frac{x^{2p+4}}{4(p+1)^2(2p+4)} + \frac{c x^2}{2p+2} - \frac{c^2}{2p x^{2p}}+C}\,,
$$
where $C$ is some other constant of integration. $u_2$ cannot be positive unless $c=0$. Hence
$$
u_1 (x) = \frac{2p+2}{x}
$$
and
$$
u_2(x) = \frac{x^{2p+1} \left ( \frac{x}{2p+2} \right )^2}{\frac{x^{2p+4}}{4(p+1)^2(2p+4)} +C}\,.
$$
etc.

Turning now to the case $u_0 \ne 0$, we have
$$
u_0 (x) = \frac{x^{-2p-1}}{c-\frac{x^{-2p}}{2 p}}
$$
for some constant of integration. Without going into details, let us simply say that $u_1$ cannot be of one sign unless $c=0$. By iterating
this argument, we deduce the expansion of $U_+$. 

Next, we examine the sequence of diffusions associated with these expansions. With the help of Table \ref{diffusionTable}, we see that the
diffusions associated with $U_-$ are:
\begin{multline}
\notag
X = \text{BES} \left (p \right) \xrightarrow[]{\text{$\phi_-(\cdot,0)$-transform}}  \text{BES} \left (p \right)  \xrightarrow[]{\text{Krein duality}} {\mathscr T}_{h_0} \left ( X \right ) 
= \text{BES} \left (-p-1 \right) \\ \xrightarrow[]{\text{$\phi_-(\cdot,0)$-transform}} \text{BES} \left (p+1 \right) \xrightarrow[]{\text{Krein duality}}  {\mathscr T}_{h_1} \circ {\mathscr T}_{h_0} \left ( X \right ) = \text{BES} \left (-p-2 \right) \;\;\text{etc}.
\end{multline}

When applied to $U_+$, using $u_{0,+} = -2p/x$, the first iteration of the algorithm yields
$$
X = \text{BES} \left (p \right) \xrightarrow[]{\text{$\phi_+(\cdot,0)$-transform}}  \text{BES} \left (-p \right)  \xrightarrow[]{\text{Krein duality}} {\mathscr T}_{h_0} \left ( X \right ) =\text{BES} \left (p-1 \right)\,.
$$
This is precisely the Ciesielski--Taylor pairing generalised by Biane \cite{Bi}.

Other well-studied diffusions that lead to simple continued fraction expansions are exponential Brownian motions and squared Bessel
processes.

\section{Diffusion in a Brownian environment with positive drift}
\label{randomSection}
Now, let $a \equiv 1$ and suppose that 
$$
W(x) = \mu x + B_x\,,\;\; \mu > 0\,,
$$
where $B$ is standard Brownian motion.
For every realisation of $W$, the process with generator (\ref{randomGenerator}) is a linear diffusion and
the algorithm of \S \ref{expansionSection} produces a continued fraction whose coefficients $u_n(x)$ are random variables.
We shall be interested in the stationary distributions of $U$ and the $u_n$.

In order to study these stationary distributions,
we require a stationary Fokker--Planck (Forward Kolmogorov) equation 
for the ``master'' system of stochastic equations
\begin{equation}
\frac{\d {\mathbf y}}{\d x} = {\mathbf  a} ( {\mathbf y}) +  B' \,{\mathbf b} {\mathbf y} \,,
\label{langevinMaster}
\end{equation}
where ${\mathbf y} :\,{\mathbb R} \rightarrow {\mathbb R}^d$, ${\mathbf a} :\,{\mathbb R}^{d} \rightarrow {\mathbb R}^d$ is a vector-valued
function with components $a_i$, and ${\mathbf b}$
is a (fixed) diagonal $d \times d$ matrix given by
$$
{\mathbf b} := \text{diag} (b_i)\,.
$$

The flow associated with Equation (\ref{langevinMaster}), interpreted in the sense
of Stratonovich \cite{RY}, defines a Feller semigroup whose
infinitesimal generator ${\mathscr G}$ is given by
\begin{equation}
{\mathscr G} f = \sum_{i=1}^d \left \{ a_i + \frac{1}{2} b_i^2 y_i  \right \} \frac{\partial f}{\partial y_i} 
+ \frac{1}{2} \sum_{i=1}^d b_i y_i \sum_{j=1}^d b_j y_j \frac{\partial^2 f}{\partial y_i \partial y_j}   \,.
\label{generator}
\end{equation}
Its adjoint ${\mathscr G}^\dag$ is the Fokker--Planck operator
\begin{equation}
{\mathscr G}^\dag f = \sum_{i=1}^d \frac{\partial}{\partial y_i} \left \{ -a_i f
+ \frac{1}{2} b_i y_i \sum_{j=1}^d \frac{\partial}{\partial y_j} \left [ b_j y_j f  \right ] \right \}\,.
\label{fokkerPlanckOperator}
\end{equation} 
If $f$ is a normalisable solution of the Fokker--Planck equation
\begin{equation}
{\mathscr G}^{\dag} f = 0\,,
\label{fokkerPlanckEquation}
\end{equation}
then it is the density of a stationary solution ${\mathbf y}$ of the stochastic equation (\ref{langevinMaster}).

\subsection{The stationary distribution of the continued fraction}
\label{continuedFractionSubsection}
The equation satisfied by the Riccati variable is
\begin{equation}
\frac{\d U}{\d x} = 2 \lambda - 2 \mu U - U^2 - 2 B' U\,.
\label{equationForU}
\end{equation}
This is of the form (\ref{langevinMaster}) with $d=1$,
$$
a = 2 \lambda - 2 \mu y - y^2\,,\quad b = -2\,.
$$
Denote by $f_U$ the density of the distribution of $U$. The Fokker--Planck equation is then
$$
\frac{\d}{\d y} \left \{ - \left ( 2 \lambda - 2 \mu y -y^2\right ) f_U + 2 y \frac{\d}{\d y} \left ( y f_U \right ) \right \} = 0\,.
$$
This has one normalisable solution, namely
\begin{equation}
f_U (y) = c y^{-\mu-1} \exp \left [ - \frac{y}{2} - \frac{\lambda}{y} \right ]\,, \;\; y > 0\,.
\label{densityOfU}
\end{equation}
This is the density of the inverse Gaussian distribution.

\subsection{Proof of Theorem \ref{randomTheorem}}
\label{randomTheoremSubsection}

The equation for the continued fraction coefficient $u_0$ is
$$
\frac{\d u_0}{\d x} = - 2 \mu u_0 - u_0^2 - 2 B' u_0\,.
$$
This is the homogeneous version of Equation (\ref{equationForU}); for $\mu > 0$, it has no normalisable solution. Hence $u_0=0$.
The equation for $u_1$ is then
\begin{equation}
\frac{\d u_1}{\d x} =  2 \mu u_1 - u_1^2 + 2 B' u_1\,.
\label{equationForu1}
\end{equation}
This is of the form (\ref{langevinMaster}) with $d=1$,
$$
a =  2 \mu y - y^2\,,\quad b = 2\,.
$$
Denote by $f$ the density of the distribution of $u_1$. The Fokker--Planck equation is then
\begin{equation}
\frac{\d}{\d y} \left \{  \left ( - 2 \mu y + y^2\right ) f + 2 y \frac{\d}{\d y} \left ( y f \right ) \right \} = 0\,.
\label{equationForf}
\end{equation}
This has one normalisable solution, namely
\begin{equation}
f (y) = c y^{\mu-1} \e^{-y/2}\,, \;\; y > 0\,.
\label{densityOfu1}
\end{equation}
This is the density of the gamma distribution.

More generally, for $d \in {\mathbb N}$, the stochastic equation satisfied by the first $d$ of the non-zero $u_n$ is of the form
(\ref{langevinMaster}) with
$$
a_i(y) = 2 (-1)^{i-1} y_i \left [ \mu + \sum_{k=1}^{i-1} (-1)^{k} y_k \right ] - y_i^2\;\;\text{and}\;\; b_i = 2 (-1)^{i-1}\,. 
$$

Equation (\ref{langevinMaster}) for the density, say $f_d$,
of the joint stationary distribution of the random variables $u_1 ,\, \ldots,\,u_{d}$ becomes
\begin{equation}
\sum_{i=1}^d \frac{\partial}{\partial y_i} \left \{ -a_i f_d + 2 (-1)^{i-1} y_i \sum_{j=1}^d (-1)^{j-1}\frac{\partial}{\partial y_j} \left [ y_j f_d \right ] \right \} 
= 0\,.
\label{equationForfd}
\end{equation}

We are now in a position to prove Theorem \ref{randomTheorem}.
Let $f$ be defined as in Equation (\ref{densityOfu1}). We will show that
the probability density function
$$
f_d (y_1,\,\ldots,\,y_d) := \prod_{i=1}^d f(y_i)
$$
solves the Fokker--Planck equation (\ref{equationForfd}).

\begin{proof}
We proceed by induction on $d$. The case $d=1$ has already been dealt with. We make the induction hypothesis: namely,
we set
$$
f_d (y_1,\,\ldots,\,y_d) := \prod_{i=1}^d f(y_i)
$$
and suppose that Equation (\ref{equationForfd}) holds for some $d$. For $d+1$,
the left-hand side of the Fokker--Planck equation may be written as the sum
$$
{\mathbf A}+{\mathbf B}+{\mathbf C}
$$ 
where
\begin{equation}
{\mathbf A} := \sum_{i=1}^d \frac{\partial}{\partial y_i} \left \{-a_i f_{d+1} + 2 (-1)^{i-1} y_i \sum_{j=1}^d (-1)^{j-1} \frac{\partial}{\partial y_j} \left [ y_j f_{d+1} \right ]
\right \},
\label{Aterm}
\end{equation}
\begin{equation}
{\mathbf B} := \sum_{i=1}^d \frac{\partial}{\partial y_i} \left \{ 2 (-1)^{i-1} y_i (-1)^{d} \frac{\partial}{\partial y_{d+1}} \left [ y_{d+1} f_{d+1} \right ]
\right \}
\label{Bterm}
\end{equation}
and
\begin{equation}
{\mathbf C} :=  \frac{\partial}{\partial y_{d+1}} \left \{-a_{d+1} f_{d+1} + 2 (-1)^{d} y_{d+1} \sum_{j=1}^{d+1} (-1)^{j-1} \frac{\partial}{\partial y_j} \left [ y_j f_{d+1} \right ]
\right \}\,.
\label{Cterm}
\end{equation}
First, we note that
$$
f_{d+1} (y_1,\,\ldots,\,y_{d+1}) = f_d(y_1,\,\ldots,\,y_d) f(y_{d+1})\,.
$$
So
\begin{multline}
\notag
{\mathbf A} = f(y_{d+1}) \sum_{i=1}^d \frac{\partial}{\partial y_i} \left \{-a_i f_{d} + 2 (-1)^{i-1} y_i \sum_{j=1}^d (-1)^{j-1} \frac{\partial}{\partial y_j} \left [ y_j f_{d} \right ]
\right \} \\
= f(y_{d+1}) \cdot 0 = 0
\end{multline}
by the induction hypothesis. Next, we remark that
\begin{multline}
\notag
{\mathbf B} = \frac{\partial}{\partial y_{d+1}} 
\sum_{i=1}^d \frac{\partial}{\partial y_i} \left \{ 2 (-1)^{i-1} y_i (-1)^{d} \left [ y_{d+1} f_{d+1} \right ]
\right \} \\
= \frac{\partial}{\partial y_{d+1}} \left \{ (-1)^d y_{d+1} f(y_{d+1}) \sum_{i=1}^d (-1)^{i-1} 2 \frac{\partial}{\partial y_i}
\left [ y_i f_d \right ]
\right \}\,.
\end{multline}
Hence
\begin{multline}
\notag
{\mathbf A}+{\mathbf B}+{\mathbf C} = \frac{\partial}{\partial y_{d+1}} 
\left \{ (-1)^d y_{d+1} f(y_{d+1}) \sum_{i=1}^d (-1)^{i-1} 2 \frac{\partial}{\partial y_i}
\left [ y_i f_d \right ]
\right . \\
\left . -a_{d+1} f(y_{d+1}) f_d + (-1)^d y_{d+1} f(y_{d+1}) \sum_{i=1}^d (-1)^{i-1} 2 \frac{\partial}{\partial y_i} \left [ y_i f_d \right ] \right . \\
\left . + 2 (-1)^d y_{d+1} (-1)^d f_d \frac{\partial}{\partial y_{d+1}} \left [ y_{d+1} f(y_{d+1}) \right ] \right \} \\
= \frac{\partial}{\partial y_{d+1}} 
\left \{ (-1)^d y_{d+1} f(y_{d+1}) \sum_{i=1}^d (-1)^{i-1} 2 \frac{\partial}{\partial y_i}
\left [ y_i f_d \right ] \right . \\
\left .  -a_{d+1} f(y_{d+1}) f_d + 2 y_{d+1}  f_d \frac{\partial}{\partial y_{d+1}} \left [ y_{d+1} f(y_{d+1}) \right ]
\right \}\,.
\end{multline}
Now,
\begin{multline}
\notag
\frac{\partial}{\partial y_i} \left [ y_i f_d \right ] = \left ( \prod_{\substack{j=1 \\ j \ne i}}^d f(y_j)\right ) \frac{\partial}{\partial y_i} \left [ y_i f(y_i) \right ]
= \left ( \frac{1}{2 y_i} \prod_{\substack{j=1 \\ j \ne i}}^d f(y_j) \right ) 2 y_i \frac{\partial}{\partial y_i} \left [ y_i f(y_i) \right ] \\
= \left ( \frac{1}{2 y_i} \prod_{\substack{j=1 \\ j \ne i}}^d f(y_j) \right ) y_i (2 \mu-y_i) f(y_i) = \frac{1}{2} (2\mu-y_i) f_d,
\end{multline}
where we have used the fact that $f(y_i)$ solves Equation (\ref{equationForf}) with $1$ replaced by $i$.
Therefore
\begin{multline}
\notag
{\mathbf A}+{\mathbf B}+{\mathbf C} = f_d \frac{\partial}{\partial y_{d+1}} \left \{ 2 (-1)^d y_{d+1} f(y_{d+1}) \sum_{i=1}^d (-1)^{i-1} (2 \mu-y_i) 
 \right . \\
\left. -a_{d+1} f(y_{d+1}) + 2 y_{d+1} \frac{\partial}{\partial y_{d+1}}  \left [ y_{d+1} f(y_{d+1}) \right ] \right \} 
\end{multline}
By using the definition of $a_{d+1}$, it is easy to verify that
$$
2 (-1)^d y_{d+1} \sum_{i=1}^d (-1)^{i-1} (2 \mu-y_i)-a_{d+1} = y_{d+1} (y_{d+1}-2\mu)\,.
$$
Hence
\begin{multline}
\notag
{\mathbf A}+{\mathbf B}+{\mathbf C} \\
= f_d \frac{\partial}{\partial y_{d+1}} \left \{  2 y_{d+1} \frac{\partial}{\partial y_{d+1}}  \left [ y_{d+1} f(y_{d+1}) \right ]  -y_{d+1} (2 \mu-y_{d+1}) f(y_{d+1}) \right \} \\ = f_d \cdot 0 = 0 
\end{multline}
and the proof is complete.
\end{proof}

If we assume that this Fokker--Planck equation admits no more
than one smooth solution, then it follows that the $u_n$ are, in the stationary regime, independent and have the same
gamma distribution. The formula (\ref{densityOfU}) reproduces the result obtained
by Letac \& Seshadri \cite{LS}. Since all the coefficients are positive, this continued fraction yields twice the
Laplace exponent $\psi_-$. In particular, we can assert that, in a Brownian environment with positive
drift, in the stationary regime, the reciprocal of the mean duration of the excursions below the starting point is gamma-distributed. The corresponding
random measure
$\sigma_-$ was studied by Marklof {\em et al.} in \cite{MTW2}; they found that the essential spectrum is $[0,\infty)$, with an empty absolutely continuous
part.

\bibliographystyle{amsplain}

\appendix

\end{document}